\documentclass[11pt]{article}
\usepackage{pifont}

\usepackage{amsfonts}
\usepackage{latexsym}
\usepackage{amsmath}
\usepackage{amssymb}
\usepackage{color}

\newtheorem{theorem}{Theorem}[section]

\newtheorem{proposition}{Proposition}[section]

\newenvironment{proof}[1][Proof]{\noindent \textbf{#1.} }{\ \ \  $\Box$}

\newtheorem{definition}{Definition}[section]
\newtheorem{remark}{Remark}[section]

\title{BSVIEs with stochastic Lipschitz coefficients and applications in finance\thanks{This work is supported by National Natural
Science Foundation of China Grant 10771122, Natural Science
Foundation of Shandong Province of China Grant Y2006A08 and National
Basic Research Program of China (973 Program, No. 2007CB814900).}}

\date{January 19 2010}

\author{Tianxiao Wang\thanks{Corresponding author, E-mail:xiaotian2008001@gmail.com}\\ \small{School of
Mathematics, Shandong University, Jinan 250100, China}}

\begin{document}

\maketitle

\begin{abstract}
This paper is concerned with existence and uniqueness of M-solutions
of backward stochastic Volterra integral equations (BSVIEs for
short) which Lipschitz coefficients are allowed to be random, which
generalize the results in \cite{Y3}. Then a class of continuous time
dynamic coherent risk measure is derived, allowing the riskless
interest rate to be random, which is different from the case in
\cite{Y3}.
\par  $\textit{Keywords:}$ Backward stochastic Volterra integral equations, Adapted
M-solutions, Dynamic coherent risk measure, Stochastic Lipschitz
coefficients
\end{abstract}



\section{Introduction}\label{sec:intro}
The literature on both static and dynamic risk measures, has been
well developed since Artzner et al \cite{ADEH} firstly introduced
the concept of coherent risk measures, see \cite{D}, \cite{FS} for
more other detailed accounts. Recently, a class of static and
dynamic risk measures were induced via $g$-expectation and
conditional $g$-expectation respectively in \cite{R1}.
$g$-expectation was introduced by Peng \cite{P1} as particular
nonlinear expectations based on backward stochastic differential
equations (BSDEs for short), which were firstly studied by Pardoux
and Peng \cite{PP1}. One nature characteristic of the above risk
measures is time-consistency (or semi-group property), however,
time-inconsistency preference usually exists in real world, see
\cite{EL}, \cite{La}, \cite{S}. As to the this case, Yong \cite{Y3}
firstly obtained a class of continuous-time dynamic risk measures,
allowing possible time-inconsistent preference, by means of backward
stochastic Volterra integral equations (BSVIEs for short) in
\cite{Y3}.

One-dimensional BSVIEs are equations of the following type defined
on $[0,T]$,
\begin{equation}
Y(t)=\psi(t)+\int_t^Tg(t,s,Y(s),Z(t,s),Z(s,t))ds-\int_t^TZ(t,s)dW(s),
\end{equation}
where $(W_t)_{t\in [0,T]}$ is a $d$-dimensional Wiener process
defined on a probability space $(\Omega ,\mathcal{F},P)$ with
$(\mathcal{F}_t)_{t\in [0,T]} $ the natural filtration of $(W_t),$
such that $\mathcal{F}_0$ contains all $P$-null sets of
$\mathcal{F}.$ The function $g:\Omega\times\Delta^c\times R\times
R^d\times R^d\rightarrow R$ is generally called a generator of (1),
here $T$ is the terminal time, and the $R$-valued
$\mathcal{F}_{T}$-adapted process $\psi(\cdot)$ is a terminal
condition; $(g,T,\psi)$ are the parameters of (1). A solution is a
couple of processes $(Y(\cdot),Z(\cdot,\cdot))$ which have some
integrability properties, depending on the framework imposed by the
type of assumptions on $g.$ Readers interested in an in-depth
analysis of BSVIEs can see \cite{Y1}, \cite{Y2}, \cite{Y3},
\cite{L}, \cite{WZ}, \cite{A} and \cite{R2}, among others.

One of the assumptions in Yong \cite{Y3} is that $r(\cdot)$ (the
interest rate) is deterministic, otherwise, it is contradicted with
the definition of translation invariance. As well known to us, in
some circumstances, it is necessary that the interest rate is
random, hence, in the current paper, we are dedicated to study the
case of the random case by giving a general version of definition
aforementioned. After that we will show a class of dynamic coherent
risk measure by means of BSVIEs, allowing the interest rate to be
random. Before doing this, we should prove the unique solvability of
M-solution, introduced by Yong \cite{Y2}, of BSVIEs under stochastic
Lipschitz condition,
\begin{equation}
Y(t)=\psi(t)+\int_t^Tg(t,s,Y(s),Z(s,t))ds-\int_t^TZ(t,s)dW(s),
\end{equation}
which generalize the result in \cite{Y3}. In addition, we claim that
this proof is much briefer than the one in \cite{Y3}.

The paper is organized as follows. Section 2 is devoted to notation.
In Section 3, we state a result of existence and uniqueness for
BSVIEs with generators satisfying a stochastic Lipschitz condition.
In Section 4, we apply the previous result to study some problems in
dynamic risk measure.

\section{Notation}
In this paper, we define several classes of stochastic processes
which we use in the sequel. We denote $\Delta^c=\{(t,s), 0\leq t\leq
s\leq T,\}$, and $\Delta=\{(t,s), 0\leq s<t\leq T,\}$. Let
$L_{\mathcal{F}_T}^2[0,T]$ be the set of $\mathcal
{B}([0,T])\otimes\mathcal{F}_T $-measurable processes
$\psi:[0,T]\times\Omega\rightarrow R$ such that
$E\int_0^T|\psi(t)|^2dt<\infty .$ We also denote
\begin{eqnarray*}
\mathcal{H}^2[0,T]=L_{\mathbb{F}}^2[0,T]\times
L^2(0,T;L_{\mathbb{F}}^2[0,T]),
\end{eqnarray*}
where $L_{\mathbb{F}}^2[0,T]$ is the set of all adapted processes $%
Y:[0,T]\times \Omega\rightarrow R$ such that
$E\int_0^T|Y(t)|^2dt<\infty ,$ and $L^2(0,T;L_{\mathbb{F}}^2[0,T])$
is the set of all processes $%
Z:[0,T]^2\times \Omega\rightarrow R$ such that for almost all $%
t\in [0,T]$, $z(t,\cdot )\in L_{\mathbb{F}}^2[0,T]$ satisfying
\begin{eqnarray*}
E\int_0^T\int_0^T|z(t,s)|^2dsdt<\infty.
\end{eqnarray*}
Now we cite some definitions introduced in \cite{Y3} and \cite{Y2}.
\begin{definition}
A mapping $\rho :L_{\mathcal{F}_T}^2[0,T]\rightarrow
L_{\mathbb{F}}^2[0,T]$ is called a dynamic risk measure if the
following hold:

1) (Past independence) For any $\psi (\cdot ),$ $\overline{\psi }(\cdot )\in L_{\mathcal{F}%
_T}^2[0,T],$ if $\psi (s)=\overline{\psi }(s),$ a.s. $\omega \in \Omega ,$ $%
s\in [t,T],$ for some $t\in [0,T),$ then $\rho (t;\psi (\cdot ))=\rho (t;%
\overline{\psi }(\cdot )),$ a.s. $\omega \in \Omega .$

2) (Monotonicity) For any $\psi (\cdot ),$ $\overline{\psi }(\cdot )\in L_{\mathcal{F}%
_T}^2[0,T],$ if $\psi (s)\leq \overline{\psi }(s),$ a.s. $\omega \in
\Omega , $ $s\in [t,T],$ for some $t\in [0,T),$ then $\rho (s;\psi
(\cdot ))\geq \rho (s;\overline{\psi }(\cdot )),$ a.s. $\quad \omega
\in \Omega,  s\in[t,T].$
\end{definition}
\begin{definition}
A dynamic risk measure $\rho :L_{\mathcal{F}_T}^2[0,T]\rightarrow L_{\mathbb{%
F}}^2[0,T]$ is called a coherent risk measure if the following hold:
1) (Translation invariance) There exists a deterministic integrable
function $r(\cdot )$ such that for
any $\psi (\cdot )\in L_{\mathcal{F}_T}^2[0,T],$%
\begin{eqnarray*}
\rho (t;\psi (\cdot )+c)=\rho (t;\psi (\cdot
))-ce^{\int_t^Tr(s)ds},\quad \omega \in \Omega ,t\in [0,T].
\end{eqnarray*}
2) (Positive homogeneity) For $\psi (\cdot )\in
L_{\mathcal{F}_T}^2[0,T]$ and $\lambda
>0,$
\begin{eqnarray*}
\rho (t;\lambda \psi (\cdot ))=\lambda \rho (t;\psi (\cdot )),\quad
a.s.\quad \omega \in \Omega ,\quad t\in [0,T].
\end{eqnarray*}
3) (Subadditivity) For any $\psi (\cdot ),$ $\overline{\psi }(\cdot )\in L_{\mathcal{F}%
_T}^2[0,T],$
\begin{eqnarray*}
\rho (t;\psi (\cdot )+\overline{\psi }(\cdot ))\leq \rho (t;\psi
(\cdot
))+\rho (t;\overline{\psi }(\cdot )),\quad \omega \in \Omega ,%
t\in [0,T].
\end{eqnarray*}
\end{definition}
\begin{definition}
Let $S\in [0,T]$. A pair of $(Y(\cdot ),Z(\cdot ,\cdot ))\in \mathcal{H}%
^2[S,T]$ is called an adapted M-solution of BSVIE (1) on $[S,T]$ if
(1) holds in the usual It\^o's sense for almost all $t\in [S,T]$
and, in addition, the following holds:
\begin{eqnarray*}
Y(t)=E^{\mathcal{F}_S}Y(t)+\int_S^tZ(t,s)dW(s).
\end{eqnarray*}
\end{definition}
\section{The existence and uniqueness with stochastic Lipschitz coefficient}
 A class of dynamic risk measures, allowing time-inconsistency preference, were induced via BSVIEs of the
 form, $t\in
 [0,T]$,
\begin{equation}
Y(t)=\psi (t)+\int_t^Tg(t,s,Y(s),Z(s,t))ds-\int_t^TZ(t,s)dW(s).
\end{equation}
In this section we will study the unique solvability of M-solution
for (3) under more weaker assumptions, i.e., allowing the
coefficients be stochastic. In addition, the proof also seems
briefer than the one in \cite{Y3}. So we introduce the standard
assumptions as follows,

(H1) Let $g:\Delta ^c\times R\times R^{d}\times \Omega \rightarrow
R$ be $\mathcal{B}(\Delta ^c\times R\times R^{d})\otimes
\mathcal{F}_T$-measurable such that $s\rightarrow g(t,s,y,z)$ is
$\mathbb{F}$-progressively measurable for all $(t,y,z)\in
[0,T]\times R\times R^{d}$ and
\begin{equation}
E\int_0^T\left( \int_t^T|g_0(t,s)|ds\right)^2dt<\infty ,
\end{equation}
where we denote $g_0(t,s)\equiv g(t,s,0,0).$ Moreover,
\begin{eqnarray}
|g(t,s,y,z)-g(t,s,\overline{y},\overline{z})| &\leq &L_1(t,s)|y-\overline{y}%
|+L_2(t,s)|z-\overline{z}|,   \\
\forall y,\overline{y} &\in &R^m,\quad z, \overline{z}\in R^{m\times
d},  \nonumber
\end{eqnarray}
where $L_1(t,s)$ and $L_2(t,s)$ are two non-negative
$\mathcal{B}(\Delta^c)\times \mathcal{F}_{T}$-measurable processes
such that for any
\begin{eqnarray}
\int_t^TL_1^2(t,s)ds<M, \quad \left( \int_t^TL_2^q(t,s)ds\right)
^{\frac 2q} <M , \quad t\in[0,T], \nonumber
\end{eqnarray}
for some constant $M$ and $\frac 1p+\frac 1q=1,$ $1<p<2$. So we
obtain the following theorem,
\begin{theorem}
Let (H1) hold, $\psi (\cdot )\in L_{\mathcal{F}_{T}}^2[0,T],$ then
(3) admits a unique M-solution in $\mathcal{H}^2[0,T]$.
\end{theorem}
\begin{proof}
Let $\mathcal{M}^2[0,T]$ be set of elements in $\mathcal{H}^2[0,T]$
satisfying: $\forall S\in [0,T]$%
\[
Y(t)=E^{\mathcal{F}_{S}}Y(t)+\int_S^tZ(t,s)dW(s).
\]
Obviously it is a closed subset of $\mathcal{H}^2[0,T]$ (see
\cite{Y2}). Due to the following inequality,
\begin{eqnarray*}
E\int_0^Te^{\beta t}dt\int_0^t|z(t,s)|^2ds\leq E\int_0^Te^{\beta
t}|y(t)|^2dt,
\end{eqnarray*}
where $\beta$ is a constant,
$(y(\cdot),z(\cdot,\cdot))\in\mathcal{M}^2[0,T],$ we can introduce a
new equivalent norm in $\mathcal{M}^2[0,T]$ as follows,
\[
\Vert (y(\cdot ),z(\cdot ,\cdot ))\Vert _{\mathcal {M}^2[0,T]}\equiv
E\left\{ \int_0^Te^{\beta t}|y(t)|^2dt+\int_0^Te^{\beta
t}\int_t^T|z(t,s)|^2dsdt\right\} ^{\frac 12}.
\]
Let us consider,
\begin{equation}
Y(t)=\psi (t)+\int_t^Tg(t,s,y(s),z(s,t))ds-\int_t^TZ(t,s)dW(s),\quad
t\in [0,T],
\end{equation}
for any $\psi (\cdot )\in L_{\mathcal{F}_{T}}^2[0,T]$ and $(y(\cdot
),z(\cdot ,\cdot
))\in \mathcal{M}^2[0,T].$ Hence BSVIE (6) admits a unique M-solution $(Y(\cdot ),Z(\cdot ,\cdot ))\in \mathcal{M}%
^2[0,T]$, see \cite{Y2}, and we can define a map $\Theta :%
\mathcal{M}^2[0,T]\rightarrow \mathcal{M}^2[0,T]$ by
\[
\Theta (y(\cdot ),z(\cdot ,\cdot ))=(Y(\cdot ),Z(\cdot ,\cdot
)),\quad \forall (y(\cdot ),z(\cdot ,\cdot ))\in \mathcal{M}^2[0,T].
\]
Let $(\overline{y}(\cdot ),\overline{z}(\cdot ,\cdot ))\in \mathcal{M}%
^2[0,T] $ and $\Theta (\overline{y}(\cdot ),\overline{z}(\cdot ,\cdot ))=(%
\overline{Y}(\cdot ),\overline{Z}(\cdot ,\cdot )).$ Obviously we
obtain that
\begin{eqnarray*}
&&E\int_0^Te^{\beta t}|Y(t)-\overline{Y}(t)|^2dt+E\int_0^Te^{\beta
t}dt\int_t^T|Z(t,s)-\overline{Z}(t,s)|^2ds \\
&\leq &CE\int_0^Te^{\beta t}\left\{ \int_t^T|g(t,s,y(s),z(s,t))-g(t,s,%
\overline{y}(s),\overline{z}(s,t))|ds\right\} ^2dt \\
&\leq &CE\int_0^Te^{\beta t}\left\{ \int_t^TL_1(t,s)|y(s)-\overline{y}%
(s)|ds\right\} ^2dt \\
&&+CE\int_0^Te^{\beta t}\left\{ \int_t^TL_2(t,s)|z(s,t)-\overline{z}%
(s,t)|ds\right\} ^2dt \\
&\leq &CE\int_0^Te^{\beta t}\left(\int_t^TL_1^2(t,s)ds\right) \int_t^T|y(s)-\overline{y}(s)|^2dsdt \\
&&+CE\int_0^Te^{\beta t}\left(\int_t^TL_2^{q'}(t,s)ds\right)^{\frac
2 {q'}}
\left(\int_t^T|z(s,t)-\overline{z}(s,t)|^{p'}ds\right)^{\frac 2 {p'}}dt \\
&\leq &CE\int_0^T|y(s)-\overline{y}(s)|^2ds\int_0^se^{\beta
t}dt+C\left[\frac1 \beta\right]^{\frac {2-p'} {p'}} E\int_0^Tds\int_t^Te^{\beta s}|z(s,t)-\overline{z}(s,t)|^2dt \\
&\leq &\frac C\beta E\int_0^Te^{\beta s}|y(s)-\overline{y}(s)|^2ds
+C\left[\frac 1 \beta\right]^{\frac {2-p'}{p'}} E\int_0^Te^{\beta t}dt\int_0^t|z(t,s)-\overline{z}(t,s)|^2ds \\
&\leq &\left(\frac C\beta +C\left[\frac 1 \beta\right]^{\frac
{2-p'}{p'}}\right)E\int_0^Te^{\beta s}|y(s)-\overline{y}(s)|^2ds,
\end{eqnarray*}
where $1<p'<2,$ $\frac {1}{p'}+\frac {1 }{q'}=1.$ Note that in above
we use the following relation, for any  $1<p'<2,$ and $r>0,$
\begin{eqnarray}
&&\left[ \int_t^T|z(s,t)-\overline{z}(s,t)|^{p'}ds\right] ^{\frac 2
{p'}} \nonumber \\
&\leq& \left[ \int_t^Te^{-rs\frac 2{2-p'}}ds\right] ^{\frac{2-p'}%
{p'}}\int_t^Te^{rs\frac 2 {p'}}|z(s,t)-\overline{z}(s,t)|^2ds \nonumber \\
&\leq &\left[ \frac 1 r\right] ^{\frac{2-p'}{p'}}\left[ \frac{2-p'}p\right] ^{%
\frac{2-p'}{p'}}e^{-rt\frac 2 {p'}}\int_t^Te^{rs\frac 2
{p'}}|z(s,t)-\overline{z}(s,t)|^2ds,
\end{eqnarray}
Then we can choose a $\beta ,$ so that the map $\Theta $ is a
contraction, and (6) admits a unique M-solution.
\end{proof}

\section{Applications in finance}
In what follows, we define
\begin{equation}
\rho(t,\psi(\cdot))=Y(t),\quad \forall t\in[0,T],
\end{equation}
where $(Y(\cdot),Z(\cdot,\cdot))$ is the unique adapted M-solution
of (3).

 Now we cite a proposition from \cite{Y3}.
\begin{proposition}
Let us consider the following form of BSVIE, $t\in [0,T]$,
\begin{equation}
Y(t)=-\psi
(t)+\int_t^T(r_1(s)Y(s)+g(t,s,Z(s,t)))ds-\int_t^TZ(t,s)dW(s).
\end{equation}
If $r_{1}(s)$ is a bounded and deterministic function, then
$\rho(\cdot)$ defined by (8) is a dynamic coherent risk measure if
$z\mapsto g(t,s,z)$ is positively homogeneous and sub-additive.
\end{proposition}
In Proposition 4.1, due to the way of defining translation
invariance in definition 2.2, $r(\cdot)$ must be a deterministic
function. In fact, we require $t\mapsto \rho(t;\psi(\cdot)+c)$ is
$\Bbb{F}$-adapted, and if we allow $r(\cdot)$ to be an
$\Bbb{F}$-adapted process, the axiom will become controversial.
Hence if we would like to consider the random case, we should
replace the axiom to a more general form. So we have,

1')There exists a $\Bbb{F}$-adapted process $Y_{0}(t)$ such that for
any $\psi (\cdot )\in L_{\mathcal{F}_T}^2[0,T],$%
\begin{eqnarray*}
\rho (t;\psi (\cdot )+c)=\rho (t;\psi (\cdot ))-Y_{0}(t),\quad
\omega \in \Omega ,t\in [0,T].
\end{eqnarray*}
Now we give a class of dynamic coherent risk measure via certain
BSVIEs. We have,
\begin{theorem}
let us consider
\begin{equation}
Y(t)=-\psi(t)+\int_t^Tl_2(t,s)Z(s,t)+l_1(t,s)Y(s)ds+\int_t^TZ(t,s)dW(s),
\end{equation}
where $l_{i}, (i=1,2)$ are two bounded processes such that $s\mapsto
l_{i}(t,s), (i=1,2)$ are $\Bbb{F}$-adapted for almost every
$t\in[0,T]$, then $\rho(\cdot)$ is a dynamic coherent risk measure.
\end{theorem}
\begin{proof}
The result is obvious, so we omit it.
\end{proof}
\begin{remark}
In Theorem 4.1, the coefficient of $Y(s)$ is random, which
generalizes the situation in Proposition 4.1. However, as we will
assert below, in this case $g$ usually can not be a general form as
in BSVIE (9).
\end{remark}
Let us consider the following two equations,
\begin{equation}
\left\{
\begin{array}{lc}
 Y^\psi(t)=-\psi (t)+\displaystyle\int_t^T(l'(t,s)Y^\psi (s)ds \nonumber
\\
\quad\quad \quad \quad +\displaystyle\int_t^Tg(t,s,Z^\psi
(s,t)))ds-\displaystyle\int_t^TZ^\psi (t,s)dW(s),\\
Y^{\psi+c}(t)=-\psi (t)-c+\displaystyle\int_t^Tl'(t,s)Y^{\psi+c}
(s)ds\nonumber \\
\quad\quad \quad \quad+\displaystyle\int_t^Tg(t,s,Z^{\psi+c}
(s,t))ds -\displaystyle\int_t^TZ^{\psi+c} (t,s)dW(s),
\end{array}
\right.
\end{equation}
where $l'$ is a process such that $s\mapsto l'(t,s)$ is
$\Bbb{F}$-adapted for almost every $t\in[0,T]$,
$(Y^{\psi}(\cdot),Z^{\psi}(\cdot,\cdot)$ and
$(Y^{\psi+c}(\cdot),Z^{\psi+c}(\cdot,\cdot)$ are the unique
M-solution of the above two BSVIEs. We denote
$$Y'(t)=Y^{\psi+c}(t)-Y^{\psi}(t),
Z'(t,s)=Z^{\psi+c}(t,s)-Z^{\psi}(t,s),\quad t,s\in[0,T].$$ Thus we
deduce that
\begin{eqnarray}
Y'(t)&=&-c+\int_t^T(g(t,s,Z^{\psi+c} (s,t))-g(t,s,Z^\psi(s,t)))ds \nonumber\\
&&+\int_t^Tl'(t,s)Y'(s)ds+\int_t^TZ'(t,s)dW(s),
\end{eqnarray}
The definition of M-solution implies
\begin{equation}
\left\{
\begin{array}{lc}
 Y^{\psi}(t)=EY^{\psi}(t)+\displaystyle\int_0^tZ^{\psi}(t,s)dW(s),\nonumber \\
Y^{\psi+c}(t)=EY^{\psi+c}(t)+\displaystyle\int_0^tZ^{\psi+c}(t,s)dW(s),\nonumber
\end{array}
\right.
\end{equation}
so we have $$Z'(t,s)=Z^{\psi+c}(t,s)-Z^{\psi}(t,s),\quad
(t,s)\in\Delta,$$ thus we can rewrite (11) as
\begin{eqnarray}
Y'(t)&=&-c+\int_t^T(g(t,s,Z^{\psi}(s,t)+Z'(s,t))-g(t,s,Z^\psi(s,t)))ds
\nonumber \\
&&+\int_t^Tl'(t,s)Y'(s)ds+\int_t^TZ'(t,s)dW(s),
\end{eqnarray}
We assume $(Y^{\psi},Z^{\psi})$ is known, then under (H1), (12)
admits a unique M-solution $(Y',Z')$, which maybe depends on
$(Y^{\psi},Z^{\psi}).$ Obviously $Y'=Y^{\psi+c}-Y^{\psi},$ which
means means that when the total wealth $\psi$ is known to be
increased by an amount $c>0$ (if $c<0,$ it means a decrease), then
the dynamic risk will be decreased (or increase) $Y'$ which maybe
depends on the total wealth $\psi$. Next we give two special cases.

(1) When $l'(t,s)$ is independent of $\omega$, let us consider the
following backward Volterra integral equation (BVIE for short),
\begin{equation}
Y^{*}(t)=-c+\int_t^Tl'(t,s)Y^{*}(s)ds,\quad t\in [0,T].
\end{equation}
Obviously by the fixed point theorem as above, (13) admits a unique
deterministic solution when $l'(t,s)$ satisfies the assumption in
(H1). It is easy to check that $(Y^{*},0)$ is the unique M-solution
of (13) when $l'$ is a deterministic function. In fact, if $Y^{*}$
is deterministic, then $Z^{*}(t,s)=0, (t,s)\in\Delta.$ After
substituting $(Y^{*},Z^{*})$ into BSVIE (12), we obtain
$$Z^{\psi+c}(t,s)=Z^{\psi}(t,s), g(t,s,Z^{\psi+c}(s,t)=g(t,s,Z^\psi(s,t),\quad (t,s)\in\Delta,$$
thus $(Y^{*},0)$ is a adapted solution of (12), moreover, it is the
unique M-solution. So in this case, $g$ can be a general form
$g(t,s,z).$ When $l'$ is independent of $t$, then
$Y^{*}(t)=-ce^{\int^T_tr(u)du},$ and we can get the result in
\cite{Y3}.

(2) When $l'$ depends on $\omega$, if the value of $(Y',Z')$ is
independent of $(Y^{\psi},Z^{\psi}),$ $g$ usually can not be the
general form $g(t,s,z)$ as above. On the one hand, the similar
result as case (1) above no longer holds. In fact, in this case,
$Y'$ usually depends on $\omega,$ i.e., there exists a $A\subseteq
[0,T]$ satisfying $\lambda(A)>0,$ such that for $t\in A,$
$P\{\omega\in\Omega, Y'(t)\neq EY'(t)\}>0,$ where $\lambda$ is the
Lebeague measure, then
$$E\int^T_0\int_0^t|Z'(t,s)|^2dsdt=E\int^T_0|Y'(t)-EY'(t)|^2dt>0,$$
which means that there must exist a set $B\subseteq \Delta,$ such
that $\forall (t,s)\in B,$
$$P(\{Z'(t,s)\neq 0\})>0, \quad \lambda(B)>0. $$ Then for a general form of $g(t,s,z),$
$$g(t,s,Z'(s,t)+Z^{\psi}(s,t))\neq g(t,s,Z^{\psi}(s,t)),$$ which means
$Y'$ maybe depends on $\psi.$ For example, if we let $l(s)=\sin
W(s),$ $c\neq0,$ and let's consider the equation below,
\begin{equation}
Y^c(t)=-c+\int_t^T\sin W(s)Y^c(s)ds+\int_t^TZ^c(t,s)dW(s).
\end{equation}
By theorem 3.1 it admits a unique M-solution. If $Y^c=EY^c,$ a.e.,
a.s., then
\begin{eqnarray}
Y^c(t) &=&-c+\int_t^TE\sin W(s)Y^c(s)ds
\end{eqnarray}
Since $E\sin W(t)=0$ implies $Y^c(t)=-c,$ then $Z(t,s)=0,$ $(t,s)\in
[0,T]^2,$ thus we have for almost any $t\in[0,T],$ $\int_t^T\sin
W(s)ds=0,$ which means that for almost any $t\in[0,T],$ $\sin
W(t)=0,$ obviously it is a contradiction. On the other hand, if
$g(t,s,Z'(s,t)+Z^{\psi}(s,t))-g(t,s,Z^{\psi}(s,t))$ is independent
of $Z^{\psi}(s,t),$ roughly speaking, the only good case for this is
$g$ is a linear function of $z$, otherwise the result will not hold.
For example, if we let $g(t,s,z)=l(t,s)z^2,$ which is a convex
function for $z$, then
$$g(t,s,Z'(s,t)+Z^{\psi}(s,t))-g(t,s,Z^{\psi}(s,t))=2l(t,s)Z'(s,t)Z^{\psi}(s,t)+Z'^{2}(s,t).$$
Obviously the value of $Y_0$ depends on $\psi$.


\begin{thebibliography}{00}
\bibitem{A}
    \newblock  A. Aman, M. N'Zi,
    \newblock   \emph{Backward stochastic nonlinear Volterra integral equations with local Lipschitz drift,}
    \newblock  Probab. Math. Stat. \textbf{25} (2005) 105--127.

\bibitem{ADEH}
    \newblock P. Artzner, F. Delbaen, J.M. Eber, D. Heath,
    \newblock \emph{Coherent measures of risk,}
    \newblock Math. Finance. \textbf{4} (1999) 203--228.

\bibitem{D}
    \newblock F. Delbaen,
         \newblock \emph{Coherent risk measures on general probability spaces,}
         \newblock  In: Sandmann, K. Schonbucher, P.J. (Eds.),
         Advances in Finance and Stochastic, Springer-Verlag,  2002, 1-37.

\bibitem{EL}
      \newblock  I. Ekeland, A. Lazrak,
      \newblock  \emph{Non-commitment in continuous time,}
      \newblock  in press.


\bibitem{FS}
     \newblock H. F$\ddot{o}$llmer, A. Schied,
     \newblock  \emph{Convex measures of risk and trading constraints.}
     \newblock  Finance. Stoch. \textbf{6} (2002) 429-447.

\bibitem{La}
     \newblock D. Laibson,
     \newblock \emph{Golden eggs and hyperbolic discounting,}
     \newblock  Q. J. Econ. \textbf{42} (1997) 443-477.

\bibitem{L}
     \newblock  J. Lin,
     \newblock \emph{Adapted solution of backward stochastic
nonlinear Volterra integral equation,}
     \newblock   Stoch. Anal. Appl. \textbf{20} (2002) 165--183.

\bibitem{P1}
    \newblock S. Peng,
    \newblock  \emph{Backward stochastic differential equations
and related $g$-expectation,} in: N. El Karoui, L. Mazliak (Eds.),
Backward Stochastic Differential Equations, in: Pitman Res. Notes
Math. Ser., Vol. 364, 1997, pp. 141--159.



\bibitem{PP1}
     \newblock E. Pardoux, S. Peng,
     \newblock \emph{Adapted solution of a backward
stochastic differential equation,}
     \newblock   Syst. Control. Lett \textbf{14} (1990) 55--61.

\bibitem{R1}
   \newblock  E. Rosazza,
   \newblock  \emph{Risk measures via g-expectation,}
   \newblock  Insurance Mathematics and Economics, \textbf{39}
   (2006) 19--34.

\bibitem{R2}
      \newblock Y. Ren,
      \newblock \emph{On solutions of Backward stochastic Volterra
      integral equations with jumps in hilbert spaces,}
      \newblock J Optim Theory Appl, in press.

\bibitem{S}
     \newblock R. H. Strotz,
     \newblock \emph{Myopia and inconsistency in dynamic utility maximization,}
     \newblock  Rev. Econ. Stud. \textbf{23} (1956) 165-180.

\bibitem{WZ}
      \newblock Z. Wang, X. Zhang,
      \newblock \emph{ Non-Lipschitz backward
stochastic volterra type equations with jumps,}
       \newblock  Stoch.Dyn. \textbf{7} (2007) 479-496.




\bibitem{Y1}
     \newblock J. Yong,
     \newblock \emph{Backward stochastic Volterra integral
equations and some related problems,},
     \newblock Stochastic Proc. Appl. \textbf{116} (2006) 779--795.

\bibitem{Y3}
     \newblock J. Yong,
     \newblock \emph{Continuous-time dynamic risk measures by backward stochastic Volterra integral
equations,},
     \newblock Appl. Anal.  \textbf{86} (2007) 1429--1442.

\bibitem{Y2}
     \newblock J. Yong,
     \newblock \emph{Well-posedness and regularity of
backward stochastic Volterra integral equation,},
     \newblock Probab. Theory
Relat. Fields. \textbf{142} (2008) 21-77.

\end{thebibliography}
\end{document}